\documentclass{amsart} 

\usepackage{amsmath,amsthm}     
\usepackage{amssymb}            
\usepackage{euscript}           
\usepackage{enumerate,calc}
\usepackage[matrix,arrow,curve,frame]{xy}    

\xymatrixcolsep{1.9pc}                          
\xymatrixrowsep{1.9pc}
\newdir{ >}{{}*!/-5pt/\dir{>}}                  

\newtheorem{thm}[subsection]{Theorem}
\newtheorem{defn}[subsection]{Definition}
\newtheorem{prop}[subsection]{Proposition}

\newtheorem{cor}[subsection]{Corollary}
\newtheorem{lemma}[subsection]{Lemma}

\theoremstyle{definition}  

\newtheorem{remark}[subsection]{Remark}

  {\end{list}}

\newcommand{\dfn}{\textbf} 

\newcommand{\mdfn}[1]{\dfn{\mathversion{bold}#1}} 

\newcommand{\cat}{\EuScript}    
\newcommand{\cC}{{\cat C}}

\newcommand{\cU}{{\cat U}}
\newcommand{\cV}{{\cat V}}
\newcommand{\cW}{{\cat W}}

\newcommand{\sSet}{s{\cat Set}}

\newcommand{\sPre}{sPre}
\newcommand{\sSh}{sSh}

\newcommand{\field}[1]  {\mathbb #1} 
\newcommand{\A}         {\field A}

\newcommand{\Z}         {\field Z}

\renewcommand{\P}         {\field P}

\DeclareMathOperator*{\colim}{colim}
\DeclareMathOperator*{\hocolim}{hocolim}

\DeclareMathOperator{\Hom}{Hom}
\DeclareMathOperator{\HOM}{\underline{Hom}}

\DeclareMathOperator{\Map}{Map}

\newcommand{\ra}{\rightarrow}                   

\newcommand{\bd}[1]{\partial\Delta^{#1}}

\newcommand{\Sm}{Sm}

\newcommand{\rea}[1]{|{#1}|}             
\newcommand{\map}{\rightarrow}

\newcommand{\ceck}[1]{\Cech(#1)}         
\newcommand{\oceck}[1]{\Cech^{o}(#1)}    
\newcommand{\oreal}[1]{\rea{\oceck{U}}}  
\newcommand{\creal}[1]{\rea{\ceck{U}}}   

\newcommand{\Cech}{\check{C}}

\numberwithin{equation}{subsection}

\begin{document}

\title{Flasque model structures for presheaves}

\author{Daniel C. Isaksen}

\address{Department of Mathematics\\ Wayne State University\\
Detroit, MI 48202}

\email{isaksen@math.wayne.edu}

\begin{abstract}
By now it is well known that there are two useful (objectwise or local)
families of
model structures on presheaves: the injective and projective.  In fact,
there is at least one more: the flasque.
For some purposes, both the projective and
the injective structure run into technical and annoying (but surmountable)
difficulties for different reasons.  The flasque model structure, which
possesses a combination of the convenient properties of both structures,
sometimes avoids these difficulties.
\end{abstract}

\maketitle

\section{Introduction}

It is well known that there are two useful classes of model structures
on simplicial presheaves: the projective and injective.
This paper is about a third class of intermediate model structures,
called the flasque structures.
Flasque presheaves arise on a regular basis in homotopical sheaf
theory \cite{BG} \cite{J1} \cite{DI}, 
so it is surprising that
these model structures have not been described previously, especially
since the proofs (as will be seen) are pretty simple if one takes
the right perspective.  The goal of this paper is to show that
the usefulness of flasque presheaves arises from the fact that there
are flasque model structures with convenient properties.

The point of the flasque model structures is that 
they share some of the advantageous properties of both the
projective and injective structures.  
In the injective model structures, it is very convenient
that the cofibrations are easy to describe and every object is cofibrant,
but the disadvantage is that the fibrations are not so easy
to describe.  
On the other hand,
the projective model structures have the advantage that the fibrations
are easy to describe, but the cost is that the cofibrations become more
complicated.  

The flasque model structures lie somewhere between the projective and
injective model structures.  By ``between'', we mean that the weak
equivalences are the same, while the class of projective cofibrations is
contained in the class of flasque cofibrations and the class of
flasque cofibrations is contained in the class of injective cofibrations.
Dually, this means that the class of injective fibrations is
contained in the class of flasque fibrations and the class of
flasque fibrations is contained in the class of projective cofibrations.

One of the confusing aspects of motivic homotopy theory is that there
is a wide choice of foundational approaches \cite{B} \cite{DHI} \cite{Hu}
\cite{MV} \cite{M}, involving presheaves or sheaves, projective or injective
cofibrations, sheaves of homotopy groups or localizations with respect
to hypercovers, etc.  The positive aspect of this circumstance is that in
any situation, one can choose the foundational approach that is best
suited to the problem at hand.

The flasque model structures are yet another possible
foundational choice.  
As an example of a result that can most easily be proven with the
flasque model structures, we show in Theorem \ref{thm:mot-htpy-gp} that
motivic stable homotopy groups
commute with filtered colimits (see also \cite{DI}). 
In fact, this problem originally motivated this project.
The key ingredient in the proof of this theorem is
that there are explicit (and verifiable) conditions for an object
to be flasque fibrant, although this description
is slightly more complicated than for projective fibrant objects.
This is entirely unlike
the injective model structures, where
there is no explicit description of the fibrant objects.

As an another example, we consider Jardine's work on the foundations
of stable motivic homotopy theory \cite{J1}.
Here the projective model structure doesn't suffice for the simple reason
that the maps $* \map \P^1$ and $* \map \A^1/\A^1-0$ are
not projective cofibrations.  
These are the usual models for the motivic sphere $S^{2,1}$, 
and it's much more difficult to define spectra with non-cofibrant spheres.

Consequently, the approach in \cite{J1} is to use the injective model 
structure, where the two models for $S^{2,1}$ are cofibrant.  
However, many of the proofs in \cite{J1} begin 
with the observation that injective fibrant presheaves are also
flasque.  This suggests that the flasque model structure is in fact
even more suitable in this context than the
injective model structure.  Note that the maps
$* \map \P^1$ and $* \map \A^1/\A^1-0$ are still flasque cofibrations.

In short, the flasque model structures should be added to the
collection of frameworks for motivic homotopy theory (along with
the projective and injective model structures) because it
is just the right tool in some situations.

This paper also answers some questions that arose in \cite{J3} about
model structures that are intermediate between the injective and
projective.  Namely, we show by example that reasonable intermediate model
structures do exist, and they are well-behaved in the sense that they
are cofibrantly generated (in fact, cellular).
Our approach to the local flasque model structure
is to perform a left Bousfield localization of the objectwise
flasque model structure, rather than to work directly with sheaves
of homotopy groups.  This localization approach gives the same
class of weak equivalences, but it is easier to verify abstract
model theoretic properties of the local model structure.

\subsection{Related work}

We mention a few papers that were absolutely essential in the development
of this project.  First, we acknowledge \cite{J2} as the seminal paper
on the subject of model structures for simplicial presheaves.

Brown and Gersten \cite{BG} defined a flasque model
structure for simplicial sheaves rather than presheaves.  In a sense,
we're doing the presheaf analogue of what they did thirty years ago.
In Theorem \ref{thm:flasque-sheaf}, we generalize part of 
their main result on the existence of certain model structures.  
Our approach has a significant advantage over \cite{BG}.  Namely, there
is no need for a Noetherian hypothesis on the Grothendieck topology.
However, without some kind of such hypothesis, it is not possible
to obtain the elegant description
of fibrant objects in terms of certain squares being homotopy cartesian.
The work of \cite{V} should also be mentioned.  This
is an axiomatization of the approach of Brown and Gersten.
Again, we get the existence of model structures in greater generality
but do not obtain the simple description of fibrant objects.

L\'arusson \cite{L} has recently used a flasque model structure in the
context of complex analysis.  His work inspired ours.  
Actually, his situation is a slight variation on the one that we consider
(see Remark \ref{rem:variation}).

The paper \cite{J3} is about model structures intermediate between
the projective and injective.
The flasque model structures are a perfect example of the kind
of structures studied there.

We should also mention that
the definition of flasque presheaf is lifted directly from \cite{J1},
which in turn is inspired by \cite{BG}.

Finally, we thank Dan Dugger for many helpful observations.

\section{Simplicial presheaves}

Recall that a simplicial presheaf on a small site $\cC$ is just a 
contravariant functor from $\cC$ to the category $\sSet$ of simplicial
sets.  We'll denote the category of presheaves by $\sPre(\cC)$.

Many of our results will be stated for arbitrary small sites, but
we're really interested in the Nisnevich
site of smooth schemes over a fixed ground scheme $S$
because that's the situation that applies to motivic homotopy theory.

The category $\sPre(\cC)$ is enriched over simplicial sets.  This
means that for every $F$ in $\sPre(\cC)$ and every simplicial set $K$,
there is a tensor $F \otimes K$ and a cotensor $F^K$ with certain
adjointness properties.  Note that $F \otimes K$ is constructed by
taking the objectwise product with $K$, and $F^K$ is constructed by 
taking the objectwise simplicial mapping space out of $K$.
The enrichment over simplicial sets also means that there are 
simplicial mapping spaces $\Map(F,G)$ for every $F$ and $G$ in
$\sPre(\cC)$.  

Every object $X$ of $\cC$ represents a simplicial presheaf of dimension
zero.  We will intentionally confuse $X$ with the presheaf that
it represents.

If $F$ is any simplicial presheaf in $\sPre(\cC)$ and $X$ is
any object of $\cC$, then $F(X)$ is naturally isomorphic to
$\Map(X, F)$.
We'll almost always use the latter
notation, but it's good to keep in mind that it's just another
name for $F(X)$.

We recollect one more construction on simplicial presheaves.
For any presheaves $F$ and $G$, there is an internal function object
$\HOM(F,G)$, which is a simplicial presheaf
such that maps $A \map \HOM(F,G)$ correspond bijectively
to maps $A \times F \map G$.

\begin{defn}
\mbox{}
\begin{enumerate}[(a)]
\item
A map $f$ in $\sPre(\cC)$ is an \mdfn{objectwise weak equivalence}
if $\Map(X,f)$ is a weak equivalence of simplicial sets for each $X$ in $\cC$.
\item
A map $f$ in $\sPre(\cC)$ is an \mdfn{injective cofibration}
if $\Map(X,f)$ is a cofibration of simplicial sets for each $X$ in $\cC$.
\item
A map $f$ in $\sPre(\cC)$ is a \mdfn{projective fibration}
if $\Map(X,f)$ is a fibration of simplicial sets for each $X$ in $\cC$.
\end{enumerate}
\end{defn}

The \dfn{projective cofibrations} are defined by a left lifting property
with respect to the maps that are both objectwise weak equivalences
and projective fibrations (i.e., objectwise acyclic fibrations).
Dually, the \dfn{injective fibrations} are defined by a right lifting
property with respect to the maps that are both objectwise weak
equivalences and injective cofibrations.

The following is a well-known result about homotopy theories
of diagram categories (see, for example, \cite{BK}, \cite{D}, and \cite{Jo}).

\begin{thm}
\label{thm:proj-inj-obj}
\mbox{}
\begin{enumerate}[(a)]
\item
The projective cofibrations, objectwise weak equivalences, and
projective fibrations form a simplicial proper cellular
model structure on $\sPre(\cC)$.
\item
The injective cofibrations, objectwise weak equivalences, and
injective fibrations form a simplicial proper cellular
model structure on $\sPre(\cC)$.
\item
The identity is a left Quillen equivalence from the projective
objectwise model structure to the injective objectwise model structure.
\end{enumerate}
\end{thm}

The proof of part (a) relies on a standard recognition principle
for cofibrantly generated categories \cite[Thm.~11.3.1]{Hi}.  The
projective generating cofibrations are maps of the form
$X \otimes \bd{n} \map X \otimes \Delta^n$, where $n \geq 0$
and $X$ belongs to $\cC$.  The projective generating acyclic
cofibrations are maps of the form
$X \otimes \Lambda^{n,k} \map X \otimes \Delta^n$, where $n \geq k \geq 0$
and $X$ belongs to $\cC$.

The proof of part (b) is quite a bit harder and relies on some complicated
set-theoretic arguments.  As a result, there is no explicit 
description of the injective fibrations.  This is the chief disadvantage
of the objectwise injective model structure.

The proof of part (c) is easy.  A projective cofibration is an injective
cofibration, and the classes of weak equivalences in the two model
categories are identical.

\section{Objectwise flasque model structure}

In this section, we work with an arbitrary small indexing category $\cC$.
For now, it is not important that $\cC$ have a Grothendieck topology.
In a few places (such as Theorem \ref{thm:obj-flasque}(c) and
Corollary \ref{cor:hom}) we need the mild hypothesis that $\cC$
contains all finite products.  In practice, this is not really
a restriction of generality.

\begin{defn}
\label{defn:union}
Let $X$ be any object of $\cC$, and let $\cU$
be a (possibly empty) finite collection of monomorphisms $U_i \map X$ 
in $\cC$.
Define $\cup \cU = \cup_{i=1}^n U_i$ to be the presheaf that
is the coequalizer of the diagram
\[
\xymatrix@1{
\coprod_{i,j} U_i \times_X U_j \ar@<0.5ex>[r]\ar@<-0.5ex>[r] &
\coprod_i U_i   }
\]
where the top arrow is projection onto the first factor and the bottom
arrow is projection onto the second factor.
\end{defn}

Beware that $\cup_{i=1}^n U_i$ is not the presheaf represented by 
the union (in $\cC$) of the $U_i$'s.  There is no guarantee that
this union even exists in $\cC$.  In fact, in our main application
to the category of smooth schemes, these unions generally don't exist:
the union of two smooth subschemes of a smooth scheme does not have
to be smooth.

Note that there is a 
canonical map $\cup \cU \map X$, and this map is a monomorphism
of presheaves.  

When $\cU$ is empty, $\cup \cU$
is equal to the empty presheaf.

Recall that if $f: F \map G$ is a map of presheaves and
$g: K \map L$ is a map of simplicial sets, then
the \mdfn{pushout product $f \Box g$} is the map
\[
G \otimes K \coprod_{F \otimes K} F \otimes L \map G \otimes L.
\]
A similar definition applies to two maps of simplicial presheaves,
except that the tensors are replaced by ordinary products.

\begin{defn}
\label{defn:flasque-gen-cofib}
\mbox{}
\begin{enumerate}[(a)]
\item
Let $I$ be the set of maps of the form $f \Box g$, where $f: \cup \cU \map X$
is induced by a finite collection of monomorphisms into an object $X$ of $\cC$
and $g$ is a generating cofibration $\bd{n} \map \Delta^n$.
\item
Let $J$ be the set of maps of the form $f \Box g$, where $f: \cup \cU \map X$
is induced by a finite collection of monomorphisms into an object $X$ of $\cC$
and $g$ is a generating acyclic cofibration $\Lambda^{n,k} \map \Delta^n$.
\end{enumerate}
\end{defn}

We will eventually show that 
$I$ serves as a set of generating flasque cofibrations,
while $J$ serves as a set of generating acyclic flasque cofibrations.

Recall that an $I$-injective is a map having the right lifting property
with respect to all elements of $I$, and an $I$-cofibration
is a map having the left lifting property with respect to all $I$-injectives
(and similarly for $J$-injectives and $J$-cofibrations).

\begin{defn}
\label{defn:flasque}
A map $f$ in $\sPre(\cC)$ is a \mdfn{flasque fibration} if 
it is a $J$-injective.
\end{defn}

By the usual adjointness arguments, a map $f: F \map G$ is flasque
if and only if the map
\[
\Map(X,F) \map \Map(X,G) \times_{\Map(\cup \cU, G)} \Map(\cup \cU, F).
\]
is a fibration of simplicial sets for all finite collections $\cU$
of monomorphisms into an object $X$ of $\cC$.

If $f$ is a flasque fibration, then $\Map(X,F) \map \Map(X,G)$ is
a fibration for all $X$ in $\cC$.  This follows from the existence of
the empty collection of monomorphisms into $X$.

A special case of Definition \ref{defn:flasque} tells us that
a presheaf $F$ is flasque fibrant if and only if the map
$\Map(X,F) \map \Map(\cup \cU, F)$ is a fibration for all finite collections
$\cU$ of monomorphisms into an object $X$ of $\cC$.

Note also that if $Y \map X$ is a monomorphism in $\cC$
and $f:F \map G$ is a flasque fibration, then the map
\[
\Map(X,F) \map \Map(X,G) \times_{\Map(Y,G)} \Map(Y, F)
\]
is a fibration.  Similarly, if $F$ is flasque fibrant, then
$\Map(X,F) \map \Map(Y,F)$ is a fibration.

\begin{remark}
\label{rem:variation}
There are many possible ways to vary the definition of flasque
fibrations.  The results proved below about flasque model structures
would hold just as well for these variations.  
We let the interested reader check that the arguments below do carry over
in the following two situations.

One possible variation is to require that a flasque fibration have the
right lifting property with respect to $\cup \cU \map X$ only for collections
$\cU$ of monomorphisms 
of size $0$ and $1$ (or for that matter, for collections up to size $n$).
Another possible variation is to only consider collections $\cU$ such that
each monomorphism $U_i \map X$ belongs to some special class of monomorphisms.
We have to assume that this special class is closed under base changes, though.

The situation of \cite[\S~12]{L} is a combination of the two
variations described above, where one considers collections
of ``Stein inclusions'' (a special kind of monomorphism) of size
$0$ or $1$.

More abstractly, in order to make the inductive proof of 
Lemma \ref{lem:acyclic-fibration} work, we need a class $C$ of finite
collections of monomorphisms such that $C$ is closed under taking
subcollections and is closed under base changes.  This last condition
means that if $\{ U_i \map X \}$ belongs to $C$ and $Y \map X$ is
any map, then $\{ U_i \times_X Y \map Y \}$ also belongs to $C$.
\end{remark}

\begin{remark}
\label{rem:traditional-flasque}
The traditional definition of a flasque presheaf $F$ of sets 
(on the Grothendieck topology of open subsets of a fixed topological space $X$)
requires that each map $F(U) \map F(V)$ be a surjection of sets for
each inclusion $V \map U$ of open sets of $X$.  When $F$ is a discrete
simplicial presheaf, our definition of a
flasque fibrant presheaf does not exactly correspond to this traditional
definition because every map between discrete simplicial sets is a fibration.
However, the philosophy is the same because intuitively fibrations
can be viewed as a certain kind of surjection, even though this is not
technically precise.
\end{remark}

\begin{defn}
A map $f$ in $\sPre(\cC)$ is a \mdfn{flasque cofibration} if it has the
left lifting property with respect to all objectwise acyclic
flasque fibrations.
\end{defn}

\begin{thm}
\label{thm:obj-flasque}
\mbox{}
\begin{enumerate}[(a)]
\item
The flasque cofibrations, objectwise weak equivalences, and
flasque fibrations form a proper cellular model structure 
on $\sPre(\cC)$.
The set of maps $I$ and $J$ (see Definition \ref{defn:flasque-gen-cofib})
serve as generating cofibrations and generating acyclic cofibrations.
\item
The identity functor on $\sPre(\cC)$ is a left Quillen equivalence
from the objectwise projective model structure to the objectwise
flasque model structure 
and from the objectwise flasque model structure to the objectwise
injective model structure.
\item If $\cC$ contains finite products, then the model structure
of part (a) is simplicial.
\end{enumerate}
\end{thm}

\begin{proof}
As for the objectwise projective model structure of 
Theorem \ref{thm:proj-inj-obj}(a), the proof is an application
of a standard recognition principle for cofibrantly generated 
model categories \cite[Thm.~11.3.1]{Hi}.  In order to apply this principle,
we have to check a few things.

First, the category $\sPre(\cC)$ contains all small limits and colimits
(they're constructed objectwise).  Second, the objectwise
weak equivalences satisfy the two-out-of-three axiom and the
retract axiom.  The rest of the hypotheses of \cite[Thm.~11.3.1]{Hi}
are proved below
in Lemmas \ref{lem:acyclic-fibration}, \ref{lem:small}, 
and \ref{lem:J-cofib}.

The simplicial structure of part (c)
(of which, as usual, only axiom SM7 is non-trivial)
is handled below by Corollary \ref{cor:SM7}.

For right properness, first note that every flasque fibration is
an objectwise fibration.  Since pullbacks are constructed objectwise
in $\sPre(\cC)$, right properness follows from right properness
of simplicial sets.

For left properness, we can apply the same argument as in the previous
paragraph if we can show that flasque cofibrations are objectwise
cofibrations.  This is proved below in Lemma \ref{lem:cofib}.
\end{proof}

The rest of this section is dedicated to proving the technical
conditions needed in the proof of Theorem \ref{thm:obj-flasque}.

\begin{lemma}
\label{lem:cofib}
A projective cofibration is a flasque cofibration, and a flasque
cofibration is an injective cofibration.  An injective fibration
is a flasque fibration, and a flasque fibration is a projective fibration.
\end{lemma}

\begin{proof}
It follows from the definitions that an objectwise acyclic flasque fibration
is an objectwise acyclic projective fibration and that
a flasque fibration is a projective fibration.
Now
projective cofibrations are determined by the
left lifting property with respect to all objectwise acyclic
projective fibrations (and similarly for flasque cofibrations),
so a projective cofibration is a flasque cofibration.

If $f: \cup \cU \map X$ is a monomorphism and 
$g: \Lambda^{n,k} \map \Delta^n$ is an acyclic cofibration,
then $f \Box g$ is an objectwise acyclic injective cofibration.
This implies that every injective fibration is a flasque fibration,
which implies in turn that every flasque cofibration is an injective
cofibration.  
\end{proof}

\begin{lemma}
\label{lem:acyclic-fibration}
A map $f: F \map G$ in $\sPre(\cC)$ 
is an objectwise acyclic flasque fibration if and only if the map
\[
\Map(X,F) \map \Map(X,G) \times_{\Map(\cup \cU, G)} \Map(\cup \cU, F)
\]
is an acyclic fibration of simplicial sets
for every finite collection $\cU$
of monomorphisms into an object $X$ of $\cC$.
Equivalently, 
a map is an objectwise acyclic flasque fibration if and only if 
it is an $I$-injective.
\end{lemma}

\begin{proof}
The second statement follows from the first by the usual adjointness
tricks.
To simplify the notation, let $\Map(\cU, f)$ be the map under
consideration.

First suppose that $\Map(\cU, f)$ is an acyclic fibration for every
finite collection $\cU$.  This immediately implies that $f$
is a flasque fibration.
Now $\Map(X,F) \map \Map(X,G)$
is an acyclic fibration for every $X$ in $\cC$ (consider the empty
collection of monomorphisms into $X$).  Thus, $f$ is an objectwise
weak equivalence.  

For the other direction, suppose that $f$ is an objectwise acyclic
flasque fibration.  From the definitions, 
this immediately implies that $\Map(\cU, f)$ is a flasque fibration.
We use induction to show that $\Map(\cU, f)$ is an objectwise weak equivalence.

When $n = 0$, the presheaf $\cup \cU$ is empty.  In this case,
$\Map(\cU, f)$ is just the map $\Map(X,F) \map \Map(X,G)$, which was
assumed to be a weak equivalence.

Now assume that the lemma has been proved 
for collections of monomorphisms of size at most $n-1$.
Let $\cU' = \{ U_1, \ldots, U_{n-1} \}$, and let
$\cU = \{U_1, \ldots, U_{n-1}, V \}$.  First note that $\cup \cU$
equals 
$\cup \cU' \coprod_{\cup \cU''} V$,
where $\cU'' = \{ U_i \times_X V \}$ is viewed as a collection of 
monomorphisms into $V$.  Thus 
$\Map(\cup \cU, F)$ equals 
\[
\Map(\cup \cU', F) \times_{\Map(\cup \cU'', F)} \Map(V, F)
\]
(and similarly for $\Map(\cup \cU, G)$).
Hence 
we have a diagram
\[
\xymatrix{
\Map(X, F) \ar[d] \\
\Map(X,G) \times_{\Map(\cup \cU, G)} \Map(\cup \cU, F) \ar[r] \ar[d] &
\Map(V, F) \ar[d] \\
\Map(X,G) \times_{\Map(\cup \cU', G)} \Map(\cup \cU', F) \ar[r]
& \Map(V,G) \times_{\Map(\cup \cU'', G)} \Map(\cup \cU'', F)       }
\]
in which the square is a pullback square.
The induction assumption tells us that the right vertical arrow
is an acyclic fibration, so the lower left vertical arrow is also.
The induction assumption also tells us that the composition of the
left column is an acyclic fibration, so the two-out-of-three axiom
allows us to conclude that $\Map(\cU, f)$ is a weak equivalence.
\end{proof}

\begin{lemma}
\label{lem:small}
The domains of the maps in $I$ or $J$ (see 
Definition \ref{defn:flasque-gen-cofib})
are $\omega$-small.
\end{lemma}

\begin{proof}
Let $F$ be the presheaf
$X \otimes K \coprod_{\cup \cU \otimes K} \cup \cU \otimes L$,
where $\cup \cU \map X$ is induced by a finite collection of monomorphisms
into an object $X$ of $\cC$ and $K \map L$ is an inclusion of finite
simplicial sets.  
It suffices to show
that each of the three presheaves $X \otimes K$, $\cup \cU \otimes K$,
and $\cup \cU \otimes L$ are $\omega$-small.  

Now $- \otimes K$
is left adjoint to the cotensor $(-)^K$, and $(-)^K$ commutes
with filtered colimits because $K$ is finite.  Thus it suffices
simply to show that $X$ and $\cup \cU$ are $\omega$-small.  Note that
$\Hom(X,F)$ is equal to the set of $0$-simplices of the simplicial
set $F(X)$; this shows that $X$ is $\omega$-small.

For $\cup \cU$, the proof is by induction on the size of $\cU$.  When
$n = 0$, the empty presheaf is certainly $\omega$-small.

Now assume that the lemma has been proved 
for collections of monomorphisms of size at most $n-1$.
Let $\cU' = \{ U_1, \ldots, U_{n-1} \}$, and let
$\cU = \{U_1, \ldots, U_{n-1}, V \}$.  First note that $\cup \cU$
equals 
$\cup \cU' \coprod_{\cup \cU''} V$,
where $\cU'' = \{ U_i \times_X V \}$ is viewed as a collection of 
monomorphisms into $V$.
As above, it suffices to observe that $\cup \cU'$, $\cup \cU''$,
and $V$ are all $\omega$-small, which follows from the induction assumption.
\end{proof}

\begin{lemma}
\label{lem:J-cofib}
If a map is a $J$-cofibration, then it is an objectwise acyclic
flasque cofibration.
\end{lemma}

\begin{proof}
Suppose that $f$ is a $J$-cofibration.
It has the left lifting property
with respect to all flasque fibrations by definition.  Therefore,
$f$ has the left lifting property with respect to all
objectwise acyclic flasque fibrations, which makes it a flasque
cofibration.  

Now we just have to show that a $J$-cofibration is an objectwise 
weak equivalence.  There are two ways of proceeding.  One way
is to just observe that $J$-cofibrations are objectwise acyclic
injective cofibrations because injective fibrations are $J$-injectives
by Lemma \ref{lem:cofib}.  The problem with this
approach is that it relies on the existence of the injective model
structure, which has a complicated proof.  We'll take a more elementary
approach.

Recall that a relative $J$-cell complex is a transfinite composition
of cobase changes of maps in $J$.  Standard arguments imply that
every relative $J$-cell complex is a $J$-cofibration.  We will first
show that relative $J$-cell complexes are objectwise weak equivalences.

By direct inspection 
a map in $J$ is an objectwise weak equivalence and an objectwise cofibration.  
Next, a cobase change of a map in $J$ is also an objectwise weak equivalence
and an objectwise cofibration.
Finally, a transfinite composition of maps that are both objectwise weak
equivalences and objectwise cofibrations
is again an objectwise weak equivalence and objectwise cofibration.

Now that we know that relative $J$-cell complexes are objectwise
weak equivalences, we consider an arbitrary $J$-cofibration $f$.
First use the small object argument (with respect to $J$)
to factor $f$ as $pi$, where 
$i$ is a relative $J$-cell complex and $p$ is a $J$-injective.
Now $f$ lifts with respect to $p$ by assumption, which implies
by the retract argument that $f$ is a retract of $i$.  Since
$i$ is an objectwise weak equivalence by the previous paragraph,
so is $f$.
\end{proof}

We used in an essential way in the
above proof that the map $\cup \cU \map X$ is a monomorphism.  Without
that, we wouldn't be able to conclude that a cobase change of a map
in $J$ is an objectwise weak equivalence.  Thus, it's critical that
we consider finite collections $\cU$ of {\em monomorphisms} into
an object $X$ of $\cC$, not all finite collections of maps
into $X$.

It's not surprising that we're forced to consider monomorphisms in $\cC$.
Otherwise, we would be producing classes of cofibrations
that aren't monomorphisms.  Although such cofibrations are not axiomatically
ruled out, in practice cofibrations are some kind of monomorphism
for most useful model categories.

\subsection{Internal Function Objects}

Our next goal is to show that the objectwise flasque model structure interacts
well with respect to internal function objects, in the same way
that the projective and injective objectwise model structures do.
We will also show that the objectwise flasque model structure is
simplicial.

\begin{lemma}
\label{lem:pushout-product}
If $K \map L$ is any cofibration of simplicial sets, then
\[
\cup \cU \otimes L \coprod_{\cup \cU \otimes K} X \otimes K \map X \otimes L
\]
is a flasque cofibration for any finite collection $\cU$ of monomorphisms
into an object $X$ of $\cC$, and it is objectwise acyclic if $K \map L$
is a weak equivalence.
\end{lemma}

\begin{proof}
The map $K \map L$ is a transfinite composition of cobase changes of
maps of the form $\bd{n} \map \Delta^n$, so the map under consideration
is a transfinite composition of cobase changes of maps of the form
$\cup \cU \otimes \Delta^n \coprod_{\cup \cU \otimes \bd{n}} X \otimes \bd{n}
   \map X \otimes \Delta^n$.  Thus the map is a relative $I$-cell complex, so
it is a flasque cofibration.

If $K \map L$ is a weak equivalence, then it is a retract of a
transfinite composition of cobase changes of maps of the form
$\Lambda^{n,k} \map \Delta^n$, so the map under consideration is a
retract of a transfinite composition of cobase changes of maps of the form
$\cup \cU \otimes \Delta^n \coprod_{\cup \cU \otimes \Lambda^{n,k}} 
   X \otimes \Lambda^{n,k} \map X \otimes \Delta^n$.  
Thus the map is a retract of a relative $J$-cell complex, so it is an
objectwise acyclic flasque cofibration.
\end{proof}

Recall that we have been working so far with an arbitrary small
indexing category $\cC$.  The next proposition requires the mild
hypothesis on $\cC$ that it contain finite products.  In applications
to motivic homotopy theory, this will be no problem.

\begin{prop}
\label{prop:pushout-product}
Let $\cC$ be an arbitrary small category that possesses finite products.
If $f: F \map G$ and $g: A \map B$ are flasque cofibrations, then
$f \Box g$ is again a flasque cofibration, and it is objectwise
acyclic if either $f$ or $g$ is.
\end{prop}

\begin{proof}
Recall that the flasque cofibrations are the retracts of relative
$I$-cell complexes,
and the objectwise acyclic flasque cofibrations 
are the retracts of relative $J$-cell complexes.
If $f$ is any map and $g$ is a retract of a transfinite composition
of cobase changes of maps in $I$, then $f \Box g$ is 
a retract of a transfinite composition
of cobase changes of maps of the form $f \Box g'$, where $g'$ 
belongs to $I$ (and similarly for $J$).  Thus it suffices to assume
that $g$ belongs to $I$ (for the first part) or to $J$ (for the second
part).  By the symmetric argument, we can also assume that $f$ belongs to $I$.

Suppose that $f$ is the map $f' \Box f''$, where
$f'$ is a map $\cup \cU \map X$ and $f''$ is the generating cofibration
$\bd{n} \map \Delta^n$.  Similarly, suppose that $g$ is $g' \Box g''$,
where $g'$ is $\cup \cV \map Y$ and $g''$ is $\bd{m} \map \Delta^m$.
Then careful inspection of the definitions shows that
$f \Box g$ is isomorphic to $h' \Box h''$, where
$h'$ is the map $f' \Box g'$ and $h''$ is the map $f'' \Box g''$.
Now $h''$ is the map
\[
\Delta^n \times \bd{m} \coprod_{\bd{n} \times \bd{m}} \bd{n} \times \Delta^m
\map \Delta^n \times \Delta^m,
\]
which is a cofibration of simplicial sets, and $h'$ is the map
$\cup \cW \map X \times Y$, where $\cW$ is the collection
consisting of
maps of the form $U_i \times Y \map X \times Y$ and maps of the form
$X \times V_i \map X \times Y$.  Thus, $h' \Box h''$ is a 
flasque cofibration by Lemma \ref{lem:pushout-product}.

The same argument works when $g$ belongs to $J$, except that $h''$
becomes
\[
\Delta^n \times \Lambda^{m,k} \coprod_{\bd{n} \times \Lambda^{m,k}} 
  \bd{n} \times \Delta^m \map \Delta^n \times \Delta^m,
\]
which is an acyclic cofibration of simplicial sets.
\end{proof}

Note that the above proof (and therefore the following two corollaries)
require us to work with arbitrary finite collections of monomorphisms
in $\cC$, not just collections of size 0 and 1.  This is the essential
reason that it is important to work with all finite collections.

\begin{cor}
\label{cor:hom}
Let $\cC$ be an arbitrary small category that possesses finite products.
If $i: A \map B$ is a flasque cofibration and $f: F \map G$ is a
flasque fibration, then the map
\[
\HOM(B, G) \map \HOM(A,G) \times_{\HOM(A,F)} \HOM(B,F)
\]
is a flasque fibration, and it is an objectwise weak equivalence
if either $i$ or $f$ is.
\end{cor}

\begin{proof}
This follows immediately by adjointness from 
Proposition \ref{prop:pushout-product}.
\end{proof}

The following corollary shows that the flasque model
structure of Theorem \ref{thm:obj-flasque} is simplicial.

\begin{cor}
\label{cor:SM7}
Let $\cC$ be an arbitrary small category that possesses finite products.
If $i: A \map B$ is a flasque cofibration and $f: F \map G$ is a
flasque fibration, then the map
\[
\Map(B, G) \map \Map(A,G) \times_{\Map(A,F)} \Map(B,F)
\]
is a fibration of simplicial sets, and it is a weak equivalence
if either $i$ or $f$ is an objectwise weak equivalence.
\end{cor}

\begin{proof}
Recall that the space of global sections of the presheaf $\HOM(F,G)$
is equal to the simplicial mapping space $\Map(F,G)$.  The
result now follows from Corollary \ref{cor:hom}.
\end{proof}

\subsection{The projective, flasque, and injective model structures are
distinct.}

We present an elementary example showing that the classes
of projective fibrations, flasque fibrations, and injective fibrations
are mutually distinct in general.

Let $\cC$ be the category with two objects $0$ and $1$ and two non-identity
morphisms from $0$ to $1$; it is indicated by the diagram
\[
\xymatrix@1{
0 \ar@<0.5ex>[r]\ar@<-0.5ex>[r] & 1.  }
\]

Now a map $F \map G$ is a projective fibration if and only if the 
maps $F_0 \map G_0$ and $F_1 \map G_1$ are fibrations.
It can be checked that $F \map G$ is a flasque fibration if and only
if the maps
$F_0 \map G_0$, $F_1 \map G_1$, and both maps
$F_1 \map F_0 \times_{G_0} G_1$ are fibrations.
Finally, $F \map G$ is an injective fibration if and only if the maps
$F_0 \map G_0$ and 
\[
F_1 \map (F_0 \times F_0) \times_{(G_0 \times G_0)} G_1
\]
are fibrations.

These three sets of conditions are distinct, as is easily verified by
elementary examples.  The reader may find it an instructive exercise
to show directly that the conditions for injective fibrations implies
the conditions for flasque fibrations.  This gives a direct
verification of Lemma \ref{lem:cofib} in this special case.

\section{Local model structures}
\label{sctn:local}

In the previous section, we worked with an index category $\cC$ that
did not necessarily have a Grothendieck topology.  From here on,
we must assume that $\cC$ does have a Grothendieck topology.

Recall that a hypercover is a certain kind of map of
simplicial presheaves 
\cite[Expos\'e V, 7.3]{SGA4} \cite{AM} \cite{DHI}.
The precise
definition is technical and not relevant for us, so we will skip it.

\begin{defn}
The \mdfn{local projective} (resp., \mdfn{flasque}, \mdfn{injective}) 
\mdfn{model structure} on $\sPre(\cC)$
is the left Bousfield localization \cite[Defn.~3.3.1]{Hi}
of the objectwise projective
(resp., flasque, injective) model structure at the class of 
all hypercovers.
\end{defn}

This is not the way that the local projective \cite{B} and local injective
\cite{J2} 
model structures are usually defined.  It is more usual (and probably more
intuitive) to define local weak equivalences with sheaves of homotopy
groups.  See \cite{DHI} for a proof that this Bousfield localization
approach gives the same class of weak equivalences for the injective
and projective model structures.  The analogous statement for the 
flasque model structure appears 
below in Theorem \ref{thm:local-wk-eq}.

It is formal that these local model structures are left proper,
cellular, and simplicial \cite[Thm.~4.1.1]{Hi}.  In fact, these model
structures are right proper.  One way to see this is to note that
the fibrations in any of these model structures are objectwise fibrations.
It follows from computing stalks that 
local weak equivalences are preserved by base
change along objectwise fibrations.

\begin{thm}
\label{thm:local-compare}
The identity functor is a left Quillen equivalence from the
local projective model structure to the local flasque model structure
and from the
local flasque model structure to the local injective model structure.
\end{thm}

\begin{proof}
In constructing the three left Bousfield localizations, one must choose
models for the homotopy colimits of the hypercovers.  It doesn't matter
which models are chosen; the localizations will be isomorphic.

Since projective cofibrations are flasque cofibrations and flasque
cofibrations are injective
cofibrations, we can choose models for the homotopy colimits in
the projective model structure, and these models will work for the
other two structures as well.  Thus, all three local model structures
are constructed by localizing at the same set of maps.

In order to show that the identity functor gives Quillen equivalences,
it suffices because of Theorem \ref{thm:obj-flasque}(b) 
to apply \cite[Thm.~3.3.20]{Hi}, which
tells us that localizations of Quillen equivalent model categories
are Quillen equivalent.
\end{proof}

The previous theorem does not guarantee that the class of 
weak equivalences in the local flasque model structure is actually
equal to the class of weak equivalences in the local
injective or local projective model structure.  However, this equality
is now not hard to prove.

\begin{thm}
\label{thm:local-wk-eq}
The local flasque weak equivalences are detected by sheaves of homotopy
groups.
\end{thm}

\begin{proof}
Let $f$ be any local flasque weak equivalence.  Factor it into a
local acyclic flasque cofibration $i$ followed by a
local acyclic flasque fibration $p$.  
Theorem \ref{thm:local-compare}
tells us that $i$ is a local acyclic injective cofibration.
Thus, $i$ induces isomorphisms on sheaves of homotopy groups.
On the other hand, $p$ is in fact an objectwise weak equivalence,
so it also induces isomorphisms on sheaves of homotopy groups.

Now suppose that $f$ is any map that induces isomorphisms on sheaves
of homotopy groups.  Factor it into a
local acyclic projective cofibration $i$ followed by a local
acyclic projective fibration $p$.  
Theorem \ref{thm:local-compare} tells us that $i$ is a local
acyclic flasque cofibration.
On the other hand, $p$ is an objectwise weak equivalence, so it
is certainly a local flasque weak equivalence.
\end{proof}

\begin{prop}
\label{prop:hom-local}
If $i: A \map B$ is a flasque cofibration and $f: F \map G$ is a
local flasque fibration, then the map
\[
\HOM(B, G) \map \HOM(A,G) \times_{\HOM(A,F)} \HOM(B,F)
\]
is a local flasque fibration, and it is a local weak equivalence
if either $i$ or $f$ is.
\end{prop}

\begin{proof}
This follows by adjointness from a local version of 
Proposition \ref{prop:pushout-product}.  The first part of that result
is no problem because the cofibrations in the local flasque model
structure are identical to the cofibrations in the objectwise
flasque model structure.

For the second part of the local version of 
Proposition \ref{prop:pushout-product}, note that
\[
B \times F \coprod_{A \times F} A \times G \map B \times G
\]
is a local weak equivalence for any pair of injective cofibrations
$A \map B$ and $F \map G$, provided at least one of them is a local weak 
equivalence.  One way to see this is to compute stalks,
noting that taking stalks commutes with finite products
and pushouts.
\end{proof}

\subsection{Simplicial sheaves}

Recall that the category $\sPre(\cC)$ of simplicial presheaves
has a full subcategory $\sSh(\cC)$ of simplicial sheaves.  The inclusion
functor has a left adjoint $a: \sPre(\cC) \map \sSh(\cC)$, which is
called the sheafification functor.  

\begin{thm}
\label{thm:flasque-sheaf}
The category $\sSh(\cC)$ has a local flasque model
structure in which the weak equivalences are the local weak
equivalences, the fibrations are the maps that are flasque fibrations
when considered as presheaves, and the cofibrations are generated
by the set $aI$ (i.e., the sheafifications of all the maps in $I$).
Sheafification is a left Quillen equivalence from the
local flasque model structure on presheaves to the local flasque
model structure on sheaves.
\end{thm}

\begin{proof}
This is an application of a general result about using adjoint functors
to translate a model structure from one category to another
\cite[Thm.~11.3.2]{Hi}.
That result has several hypothesis, only one of which is not immediately
obvious.  Namely, we must explain why relative $aJ$-cell complexes
are local weak equivalences.  This follows from the observation that
the local weak equivalences are closed under transfinite compositions.

To see that $a$ is a Quillen equivalence, one just needs to use that
the natural map $F \map aF$ is always a local weak equivalence
for any simplicial prehseaf $F$ \cite[Lem.~2.6]{J2}.
\end{proof}

The previous theorem reproves the part of \cite{BG} 
that establishes a model structure on the category of sheaves on a
Noetherian topological space.  At first glance, 
Theorem \ref{thm:flasque-sheaf} does not appear to apply because
\cite{BG} makes no mention of finite collections of monomorphisms.
However, in the Grothendieck topology of open sets of a fixed topological
space, the object $\cup \cU$ is always representable.

In fact, Theorem \ref{thm:flasque-sheaf} is much more general because there
is no Noetherian hypothesis on the site.  However, we cannot recover
the simple identification of fibrant objects in terms of certain homotopy
pullback squares.  This identification is an important part of
\cite{BG} (and its generalization in \cite{V}).

\subsection{Nisnevich local model structures}

From now on,
we specialize to the case of the Nisnevich topology on the category
$\Sm_S$ of smooth schemes over a ground scheme $S$
because this is the situation that arises in motivic homotopy theory.
In this context,
it is not necessary to localize with respect to all hypercovers
as in Section \ref{sctn:local}.
Recall that there is a certain class of elementary Nisnevich squares
\cite[Defn.~3.1.3]{MV}
\[
\xymatrix{
U \times_X V \ar[r] \ar[d] & V \ar[d] \\
U \ar[r] & X.                          }
\]
The full definition is not so important, but we will use the fact that
$U \map X$ (and thus also $U \times_X V \map V$) is a monomorphism.

The following theorem is proved in \cite[Lem.~4.2]{B}.

\begin{thm}
\label{thm:Nis-proj}
The local projective model structure on $\sPre(\Sm_S)$ is the
left Bousfield localization of the objectwise projective model
structure at the set of consisting of maps $P \map X$, where 
$P$ is the homotopy pushout of
\[
\xymatrix@1{U & U \times_X V \ar[r] \ar[l] & V}
\]
for all elementary Nisnevich squares.
\end{thm}

The homotopy pushouts in the previous theorem are a bit of an annoyance.
It would be nicer to have a more concrete description of the maps
with respect to which we are localizing.  In both the flasque
and injective cases, such a description is possible.

\begin{thm}
\label{thm:Nis-flasque}
The local flasque (resp., injective) model structure on $\sPre(\Sm_S)$ is the
left Bousfield localization of the objectwise flasque (resp., injective) model
structure at the set of maps
\[
\xymatrix@1{U \coprod_{U \times_X V} V \map X}
\]
for all elementary Nisnevich squares.
\end{thm}

\begin{proof}
Recall that localizations of Quillen equivalent model categories
are Quillen equivalent \cite[Thm.~3.3.20]{Hi}.
Therefore,
if we localize the flasque or injective objectwise model structure
at the set of maps $P \map X$ (as in the previous theorem), then 
Theorems \ref{thm:obj-flasque}(b) and \ref{thm:Nis-proj} tell us that we
get the desired local model structure.
But in the flasque and
injective model structures, $U \times_X V \map V$ is a cofibration,
so the ordinary pushout 
$U \coprod_{U \times_X V} V$ is a model for the homotopy pushout.
\end{proof}

The above theorem demonstrates 
one of the ways in which the flasque (or injective) model structure is
more convenient than the projective model structure.  The localization
is slightly more concrete --- there aren't any homotopy pushouts involved.

Theorem \ref{thm:Nis-proj} formally leads to the following characterization
\cite[Lem.~4.1]{B}.
An object $F$ of $\sPre(\Sm_S)$ is local projective 
fibrant if and only if:
\begin{enumerate}
\item
$F$ is projective fibrant;
\item
\[
\xymatrix{
F(X) \ar[r] \ar[d] & F(U) \ar[d] \\
F(V) \ar[r] & F(U \times_X V)     }
\]
is homotopy cartesian for all elementary Nisnevich squares.
\end{enumerate}

Because of Theorem \ref{thm:Nis-flasque}, 
we can improve on this characterization
for the flasque and injective model structures.

\begin{cor}
\label{cor:local-flasque-fib}
An object $F$ of $\sPre(\Sm_S)$ is local flasque (resp., injective)
fibrant if and only if:
\begin{enumerate}
\item
$F$ is flasque (resp., injective) fibrant;
\item
the map 
\[
F(X) \map F(U) \times_{F(U \times_X V)} F(V)
\]
is an acyclic fibration for all elementary Nisnevich squares.
\end{enumerate}
\end{cor}

\begin{proof}
This follows from formal properties of left Bousfield
localizations \cite[Thm.~4.1.1(2)]{Hi} (see \cite[Cor.~7.1]{DHI} for
a similar argument).
\end{proof}

The previous result shows exactly why the flasque model structure
is a good compromise between the projective and injective model structures.
Condition (1) for the injective model structure is very complicated.
However, condition (1) for
the flasque structure is much simpler.
Condition (2) doesn't work for the projective structure because
a messier homotopy pullback statement is required.

\section{Motivic model structures}

\begin{defn}
The motivic projective (resp., flasque, injective) 
model structure on $\sPre(\Sm_S)$
is the left Bousfield localization of the local projective
(resp., flasque, injective) model structure at the set of 
maps $X \map X \times \A^1$.
\end{defn}

The maps $X \map X \times \A^1$ are induced by the inclusion $0 \map \A^1$.

It is formal that these motivic model structures are left proper,
cellular, and simplicial \cite[Thm.~4.1.1]{Hi}.  In fact, they are also
right proper, but this requires some extra work \cite[2.2.7]{MV} (beware
that this reference reverses the definition of left proper and right proper).

Similarly to the local projective model structure,
an object $F$ is motivic projective 
fibrant if and only if:
\begin{enumerate}
\item
$F$ is projective fibrant;
\item
\[
\xymatrix{
F(X) \ar[r] \ar[d] & F(U) \ar[d] \\
F(V) \ar[r] & F(U \times_X V)     }
\]
is homotopy cartesian for all elementary Nisnevich squares;
\item
$F(X \times \A^1) \map F(X)$ is a  weak equivalence for all $X$ in $\cC$.
\end{enumerate}

As in Corollary \ref{cor:local-flasque-fib}, 
we can improve on this characterization
for the flasque and injective model structures.

\begin{cor}
\label{cor:mot-flasque-fib}
An object $F$ is motivic flasque (resp., injective)
fibrant if and only if:
\begin{enumerate}
\item
$F$ is flasque (resp., injective) fibrant;
\item
the map
\[
F(X) \map F(U) \times_{F(U \times_X V)} F(V)
\]
is an acyclic fibration;
\item
$F(X \times \A^1) \map F(X)$ is an acyclic fibration for all $X$ in $\cC$.
\end{enumerate}
\end{cor}

\begin{proof}
As in the proof of Corollary \ref{cor:local-flasque-fib},
this follows from \cite[Thm.~4.1.1(2)]{Hi}.
\end{proof}

%
%

We present one consequence of Corollary \ref{cor:mot-flasque-fib}, 
which will be
important later.  Note that the injective analogue of this statement
is false.

\begin{prop}
\label{prop:unstable-mot-fib-colim}
The class of motivic flasque fibrant objects is closed under
filtered colimits.
\end{prop}

\begin{proof}
Condition (1) of Corollary \ref{cor:mot-flasque-fib} is preserved by 
filtered colimits
because the generating objectwise acyclic flasque cofibrations (i.e.,
the set $J$ of Definition \ref{defn:flasque-gen-cofib})
have compact domains.  
Conditions (2) and (3) are straightforward (see also \cite[\S~1.3]{J1}).
\end{proof}

\begin{prop}
\label{prop:hom-mot}
If $A$ is flasque fibrant and $F$ is
motivic flasque fibrant, then 
$\HOM(A,F)$ is also motivic flasque fibrant.
\end{prop}

\begin{proof}
From Proposition \ref{prop:hom-local}, we know that
$\HOM(A,F)$ is local flasque fibrant.
By \cite[Thm.~4.1.1(2)]{Hi}, we only need to show that the map
$\HOM(A,F) \map *$ has the right lifting property with respect to every
map of the form
\[
X \otimes \Delta^n \coprod_{X \otimes \Lambda^{n,k}} 
  (X \times \A^1) \otimes \Lambda^{n,k}
\map (X \times \A^1) \otimes \Delta^n.
\]
Now this map is an objectwise acyclic flasque cofibration, so
the lifting property is satisfied.
\end{proof}

In analogy to Corollary \ref{cor:hom} and Proposition \ref{prop:hom-local},
it is natural to consider the map
\[
\HOM(B, G) \map \HOM(A,G) \times_{\HOM(A,F)} \HOM(B,F)
\]
for $i:A \map B$ a flasque cofibration
and $f:F \map G$ a motivic flasque fibration.  However, we do not yet
know how to show that this map is a motivic flasque fibration.

\section{Stable motivic model structures}

\subsection{Pointed model structures}

In stable motivic homotopy theory, we're actually interested in
the pointed version of what we've been discussing so far, i.e., the 
category $\sPre(\Sm_S)_*$ of presheaves of pointed simplicial sets.
Of course, this is just an undercategory of $\sPre(\Sm_S)$,
so we obtain formally the pointed objectwise flasque model structure,
the pointed local flasque model structure, and the pointed
motivic flasque model structure \cite[Thm.~7.6.5]{Hi}.

These pointed model structures
are still simplicial model categories with the obvious pointed
analogue of tensor.  Namely, if $F$ is any pointed presheaf and $K$
is any simplicial set, then $F \otimes K$ is the pointed presheaf
defined by $(F \otimes K)(U) = F(U) \wedge K_+$.

Recall that the fibrant objects in the pointed categories are
simply the pointed presheaves that are fibrant in the corresponding
unpointed category.  On the other hand, the cofibrant objects
are different; instead of studying $\phi \map A$, we now look at
$* \map A$.  We're going to need to know that some of the basic
objects of motivic homotopy theory are flasque cofibrant
pointed presheaves.

\begin{lemma}
If $X \map Y$ is any monomorphism in $\Sm_S$, then 
$Y/X$ is a flasque cofibrant pointed presheaf.
\end{lemma}

\begin{proof}
Just consider the pushout square
\[
\xymatrix{
X \ar[r] \ar[d] & Y \ar[d] \\
\mbox{}* \ar[r] & Y/X,                  }
\]
in which the top arrow is a flasque cofibration.  Therefore,
the bottom arrow is also a flasque cofibration.
\end{proof}

The objects mentioned in the following lemma are central to the
construction of stable motivic homotopy theory.

\begin{lemma}
\label{lem:cofib-sphere}
In $\sPre(\Sm_S)_*$, the following pointed presheaves are flasque cofibrant:
\begin{enumerate}
\item
The presheaf $\P^1$, pointed at $\infty$ (or at any other rational point).
\item
The presheaf $\A^1-0$, pointed at $1$ (or at any other rational point).
\item
The presheaf $\A^1/\A^1-0$.
\item
The presheaf $\Sigma (\A^1-0)$.
\end{enumerate}
\end{lemma}

\begin{proof}
The first three examples are special cases of the previous lemma.
For the last one, we just need to show that the map
\[
* \otimes S^1 \coprod_{* \otimes *} (\A^1-0) \times * \map (\A^1-0) \otimes S^1
\]
is a flasque cofibration (by the definition of suspension).
This is exactly what Lemma \ref{lem:pushout-product} says.
\end{proof}

The category $\sPre(\Sm_S)_*$ still has internal function objects,
which we denote by $\HOM_*(F,G)$.  However, we can't prove
the pointed analogues of Proposition \ref{prop:pushout-product} 
and Corollary \ref{cor:hom} because
that would require taking smash products of representable
functors instead of products of representable functors. 
Nevertheless, we have the following useful result.

\begin{prop}
\label{prop:pointed-hom}
If $A$ is a flasque cofibrant pointed presheaf and $F$ is a 
flasque (resp., local flasque, motivic flasque) fibrant pointed presheaf,
then $\HOM_*(A,F)$ is also a flasque (resp., local flasque, motivic flasque)
fibrant pointed presheaf.
\end{prop}

\begin{proof}
Consider the pullback square
\[
\xymatrix{
\HOM_*(A,F) \ar[r] \ar[d] & \HOM(A,F) \ar[d] \\
\mbox{}* \ar[r] & \HOM(*,F).                 }
\]
We are interested in the left vertical arrow, which is a base change
of the right vertical arrow.  Thus, we just have to show that the
right vertical arrow is a flasque (resp., local flasque, motivic flasque)
fibration.  For the first two cases, this follows from 
Corollary \ref{cor:hom} and Proposition \ref{prop:hom-local}.

The third (motivic) case requires only slightly more work.
First note that both objects $\HOM(A,F)$ and $\HOM(*,F)$ are
motivic flasque fibrant by Proposition \ref{prop:hom-mot}.
We already know that the map between them is a flasque fibration
from the previous paragraph.
From \cite[Prop.~3.3.16(1)]{Hi} (a general result describing
fibrations between fibrant objects in a localization), we 
conclude that the map is in fact a motivic flasque fibration.
\end{proof}

Recall that for any pointed presheaf $F$, the 
presheaf $\Omega^{2,1} F$ is defined to be $\HOM_*(T, F)$,
where $T$ is some chosen model for the motivic sphere $S^{2,1}$.
For different purposes, authors have chosen $T$ to be 
$\P^1$, $\A^1/\A^1-0$, or $\Sigma (\A^1-0)$.  In the following
corollary, it doesn't matter which of these three are used.

\begin{cor}
\label{cor:loops}
If $F$ is a motivic flasque fibrant pointed presheaf, then
$\Omega^{2,1} F$ is also a motivic flasque fibrant pointed presheaf.
\end{cor}

\begin{proof}
This is a straightforward combination of Lemma \ref{lem:cofib-sphere} with
Proposition \ref{prop:pointed-hom}.
\end{proof}

The above corollary is a great example of the advantage of the flasque
model structures over the projective model structures.  Since
$* \map \P^1$ and $* \map \A^1/(\A^1-0)$ are not projective cofibrations,
one has to choose carefully a model for $S^{2,1}$ when working with
the projective structures.  With the flasque (or injective) structures,
there is no such problem.

\subsection{Motivic spectra}

One of the key applications of the motivic flasque model structure is
for stable motivic homotopy theory.  Normally people use
the injective structure because they need $* \map \P^1$ 
(or $* \map \A^1/\A^1-0$)
to be a cofibration.
Of course, this map is also a flasque cofibration, so the general machinery
for producing stable model structures \cite{Ho} works fine.

To illustrate this point, we will describe a theorem that can be proved
only with the flasque model structure (as far as we know).  The proof of
the theorem will require certain properties that neither the
projective nor injective model structures possess.

We briefly review the construction of naive Bousfield-Friedlander
motivic spectra. 
All of what follows would work just as well
for motivic symmetric spectra, but we use
naive motivic spectra because they are easier to describe and 
they suffice for our specific purpose.  

A naive motivic spectrum $E$ is a sequence $\{ E_n \}$ of pointed
simplicial presheaves together with structure maps 
$E_n \map \Omega^{2,1} E_{n+1}$.  The definition of
maps of naive motivic spectra is obvious.

In order to construct the stable motivic model structure,
we can start with the projective, flasque, or injective model structure
on pointed simplicial presheaves.  We'll start with the flasque
model structure.

We only need a few facts about the stable motivic model structure on
naive motivic spectra.
First, a spectrum $E$ is fibrant if and only if each $E_n$ is flasque fibrant
and the maps $E_n \map \Omega^{2,1} E_{n+1}$ are (unstable) motivic
weak equivalences.
Second, a map $E \map F$ of fibrant spectra 
is a stable motivic weak equivalence if and only if each $E_n \map F_n$
is a motivic weak equivalence of pointed simplicial presheaves.
Third, the standard models for the sphere spectra $S^{p,q}$ are
cofibrant spectra.

\begin{prop}
\label{prop:stable-mot-fib-colimit}
Fibrant motivic spectra are closed under filtered colimits.
\end{prop}

\begin{proof}
Let $\{ E^i \}$ be filtered system of naive motivic spectra, and
let $E$ be $\colim_i E^i$.  Then $E_n$ equals $\colim_i E^i_n$, so
it is a motivic flasque fibrant pointed simplicial presheaf by Proposition
\ref{prop:unstable-mot-fib-colim}.  

By Corollary \ref{cor:loops}, $\Omega^{2,1} E^i_{n+1}$ is a
motivic flasque fibrant pointed simplicial presheaf, so the structure map
$E^i_n \map \Omega^{2,1} E^i_{n+1}$ is a motivic weak equivalence
between motivic flasque fibrant pointed simplicial presheaves.  Therefore,
it is in fact an objectwise weak equivalence (this is a general
property of localizations).  It follows that
$\colim_i E^i_n \map \colim_i \Omega^{2,1} E^i_{n+1}$ 
is also an objectwise weak equivlence and hence a motivic weak equivalence.
This last map is just the structure map $E_n \map \Omega^{2,1} E_{n+1}$
because $\Omega^{2,1}( -)$ commutes with filtered colimits.
\end{proof}

\begin{thm}
\label{thm:mot-htpy-gp}
Let $\{E^i\}$ be a filtered system of motivic spectra.
The natural map $\colim_i \pi_{p,q} E^i \map \pi_{p,q} (\hocolim_i E^i)$
is an isomorphism for all $p$ and $q$ in $\Z$.
\end{thm}

\begin{proof}
We may assume that each $E^i$ is a fibrant motivic spectrum.  By
Proposition \ref{prop:stable-mot-fib-colimit}, so is $\colim_i E^i$.
Consider the map
\[
\colim_i \pi_0 \Map(S^{p,q}, E^i) \map \pi_0 \Map(S^{p,q}, \colim_i E^i).
\]
This is an isomorphism because $\Map(S^{p,q}, -)$ commutes with
filtered colimits.

It only remains to show that the natural map 
$\hocolim_i E^i \map \colim_i E^i$ is a stable motivic weak equivalence.
It suffices to show that if $\{E^i\}$ and $\{F^i\}$ are
filtered systems of motivic fibrant spectra
and each $E^i \map F^i$ is a
stable motivic weak equivalence, 
then $\colim_i E^i \map \colim_i F^i$ is a stable motivic 
weak equivalence.  

Since each $E^i \map F^i$ is a stable motivic weak equivalence
between fibrant spectra, each $E^i_n \map F^i_n$ is a motivic
weak equivalence between motivic flasque fibrant pointed simplicial
presheaves.
As observed in the previous proof, each map
$E^i_n \map F^i_n$ is in fact an objectwise weak equivalence.
It follows
that $\colim_i E^i_n \ra \colim_i F^i_n$ is still an objectwise
weak equivalence, and this implies that $\colim E_i \ra
\colim F_i$ is a stable motivic weak equivalence.
\end{proof}

\bibliographystyle{amsalpha}

\end{document}